\documentclass[12pt, a4paper]{article}

\usepackage{amsmath}\multlinegap=0pt
\usepackage{amssymb}

\usepackage{epstopdf}
\usepackage[francais,english]{babel}
\providecommand\mathbb{\bf}
\newcommand\R{{\mathbb R}}

\newcommand\pref[1]{(\ref{#1})}

\let \eps\varepsilon
\newtheorem{thm}{Theorem}

\newtheorem{prop}{Proposition}
\newtheorem{lem}{Lemma}

\newtheorem{defi}{Definition}

\newcounter{Remark}
\newenvironment{Remark}
      {\par\medbreak \refstepcounter{Remark}
       \noindent\textbf{Remark~\arabic{Remark}.}}
      {\hbox{ }}
\newenvironment{proof}[1][]
      {\par\medbreak{\noindent\bfseries Proof#1.\quad}}
      {\hbox{}\hfill\fbox{\ }\bigbreak}
\newcounter{steps}

\newcommand\F{{\cal F}}

\def\<#1,#2>{\left<#1,#2\right>}

\title{Hamilton-Jacobi-Bellman equations for the optimal control of a state equation with  memory}
\author {G.~Carlier, R. Tahraoui 
\thanks{\scriptsize\ Universit\'e Paris Dauphine, CEREMADE, Pl. de Lattre de Tassigny, 75775 Paris Cedex 16, FRANCE \texttt{carlier@ceremade.dauphine.fr, tahraoui@ceremade.dauphine.fr} }}

\begin{document}

\maketitle

\begin{abstract}
This article is devoted to the optimal control of state equations with memory of the form:
\[\begin{split}
\dot{x}(t)=F\left(x(t),u(t), \int_0^{+\infty} A(s) x(t-s) ds\right), \; t>0,\\
\mbox{ with initial conditions }  x(0)=x, \; x(-s)=z(s), s>0.
 \end{split}\]
Denoting by  $y_{x,z,u}$ the solution of the previous Cauchy problem and:
\[v(x,z):=\inf_{u\in V} \left\{ \int_0^{+\infty} e^{-\lambda s } L(y_{x,z,u}(s), u(s))ds \right\}\]
where $V$ is a class of admissible controls, we prove that $v$ is the only viscosity solution of an Hamilton-Jacobi-Bellman equation of the form:
\[\lambda v(x,z)+H(x,z,\nabla_x v(x,z))+\<D_z v(x,z), \dot{z} >=0\]
in the sense of the theory of viscosity solutions in infinite-dimensions of M. Crandall and P.-L. Lions.
\end{abstract}


\clearpage

\section{Introduction}

The optimal control of dynamics with memory  is an issue that naturally arises in many different applied settings both in engineering and decision sciences. It is typically the case when studying the optimal performances of  a system in which the response to a given input occurs not instantaneously but only after a certain elapse of time. To cite some recent related contributions, in a stochastic framework,  we refer to I. Elsanosi, B. \O{}ksendal, A. Sulem \cite{eos}, for applications to mathematical finance, and to F. Gozzi and C. Marinelli  \cite{gm} for applications to advertising modelling. In the deterministic case, we refer to R. Boucekkine et al. \cite{bou}  for a generalization of Ramsey's economic growth model with memory effects and, in the field of biosciences modelling, we refer to the survey of C.T. H. Baker et al. \cite{baker}.

\smallskip

The aim of the present article is to study, by dynamic programming arguments, the optimal control of (deterministic) state equations with memory. For the sake of simplicity, we will restrict the analysis to (finite-dimensional) dynamics of the form:
\[\dot{x}(t)=F\left(x(t),u(t), \int_0^{+\infty} A(s) x(t-s) ds\right), \; t>0,\]
 with initial conditions $x(0)=x$ and  $x(-s)=z(s)$, $s>0$.
We will also focus on the discounted infinite horizon problem:
\begin{equation}\label{dvaleur}
v(x,z):=\inf_{u\in V} \left\{ \int_0^{+\infty} e^{-\lambda s } L(y_{x,z,u}(s), u(s))ds \right\}
\end{equation}
where $V$ is some admissible class of controls. 
\smallskip
Of course, there are other forms of memory effects than the one we treat here: systems with lags or with deviating arguments for instance (see for instance \cite{ct}, \cite{ts}, \cite{ts2} and the references therein). 

\smallskip

As is obvious from \pref{dvaleur}, the value function depends not only on the current state of the system $x$ but also on the whole past of the trajectory i.e. $z$ (note that we have not required $z(0)=x$ in \pref{dvaleur}). Hence the state space for problem \pref{dvaleur} is infinite-dimensional. Note that if we had imposed an additional continuity condition ensuring $x=z(0)$, then the value function would have been a function of the past $z$ only. There are several reasons why we have not adopted this point of view and have prefered to write everywhere $x$ and $z$ as if they were independent variables. The main one, is that it enables to understand the tight connections between the control problem \pref{dvaleur} and the Hamilton-Jacobi-Bellman equation:
\begin{equation}\label{hjbd}
\lambda v(x,z)+H(x,z,\nabla_x v(x,z))+\<D_z v(x,z),\dot{z}>=0 \mbox{ if } z(0)=x.
\end{equation}
The previous equation presents several difficulties. The first one is of course its infinite-dimensional nature. The second one comes from the presence of the time derivative of $z$, $\dot{z}$ in the equation and the third one from the restriction $x=z(0)$. In a series of articles (\cite{cranlio, cranlio2, cranlio3, cranlio4, cranlio5, cranlio6}), M. Crandall and P.-L. Lions developed a general theory of viscosity solutions in infinite dimensions. This theory is of course of particular interest for the optimal control of  infinite-dimensional systems. In such problems,  the Hamilton-Jacobi equation frequently contains an unbounded linear term (as in \pref{hjbd}) and in (\cite{cranlio4, cranlio5, cranlio6}),  M. Crandall and P.-L. Lions  showed how to overcome this additional difficulty. The main contribution of the present paper is to show, in a rather simple and self-contained way, how the theory of viscosity solutions in infinite dimensions of M. Crandall and P.-L. Lions can be applied to fully characterize the value function \pref{dvaleur} as the unique solution of \pref{hjbd}. For the sake of simplicity, we will work in a hilbertian framework i.e. in the state space $E:=\R^d\times L^2(\R_+,\R^d)$ and defining:
\[E_0:=\{(z(0),z), \; z\in H^1(\R_+, \R^d)\},\] 
$v$ will be said to be a viscosity subsolution of \pref{hjbd} if for every $(x_0,z_0)\in \R^d\times L^2$ and every $\phi\in C^1(\R^d\times L^2,\R)$ such that $v-\phi$ has a local maximum (in the sense of the strong topology of $\R^d\times L^2$)
 at $(x_0,z_0)$, one has:
 \[\lambda v(x_0,z_0)+H(x_0,z_0,\nabla_x \phi(x_0,z_0))+\mathrm{liminf}_{(x,z)\in E_0 \rightarrow (x_0,z_0) } \<D_z \phi(x,z),\dot{z}>\leq 0.\]
 Supersolutions of \pref{hjbd} are defined in a similar way. Now, a convenient way to study \pref{hjbd} is to rewrite it as an Hamilton-Jacobi equation with an unbounded linear term as in  M. Crandall and P.-L. Lions \cite{cranlio4, cranlio5, cranlio6}. Namely, defining $\alpha=(x,z)$, the equation reads as:
\begin{equation}\label{hjbd2}
\lambda v(\alpha) +H(\alpha, \nabla_x v(\alpha))-x\cdot \nabla_x v(\alpha)+\<T^*(\alpha), Dv(\alpha)>=0, \alpha\in D(T^*).\end{equation}
Where  $T$ is the linear unbounded operator on $E$ with domain $D(T)=\R^d\times H^1$  defined by
\begin{equation}\label{defofTb}
T(y,w):=(y-w(0),-\dot{w}),\; \forall (y,w)\in D(T).
\end{equation}
So that its adjoint, $T^*$ has domain $D(T^*)=E_0$ and is given by
\begin{equation}
T^*(x,z):=(z(0),\dot{z})=(x,\dot{z}),  \; \forall (x,z)\in D(T^*)=E_0.
\end{equation}\label{defofTadb}

\smallskip

Section \ref{sec1} is devoted to some preliminaries on the Cauchy problem and continuity properties of the value function. Section \ref{sec2} concerns the dynamic programming principle. In section \ref{sec3}, we identify the Hamilton-Jacobi-Bellman equation of the problem and establish that the value function  is a viscosity solution of this equation. In section \ref{sec4}, we prove a comparison result. Finally, in section \ref{sec5}, we end the paper by some concluding remarks.

\section{Assumptions and preliminaries}\label{sec1}

\subsection{On the Cauchy problem}

Let $K$ be a compact metric space, we define the set of admissible controls $V$ as the set of measurable functions on $(0,+\infty)$ with values in $K$. For $z\in L^2:=L^2((0,+\infty),\R^d)$, $x\in \R^d$ and $u\in V$ an admissible control, we consider the following controlled equation
\begin{equation}\label{state}
\dot{x}(t)=F\left(x(t),u(t), \int_0^{+\infty} A(s) x(t-s) ds\right), \; t>0,
\end{equation}
together with the boundary conditions:
\begin{equation}\label{bc}
x(0)=x,\; x(-s)=z(s), \; s>0.
\end{equation}
In the paper, $d$ and $k$ are given positive integers and we will always assume the following on the data $A$ and  $F$:

\begin{itemize}
\item {\bf{(H1)}} $F\in C^0(\R^d\times K\times \R^k, \R^d)$ and there exists a constant $C_1\geq 0$ such that:
\begin{equation}\label{flip}
\vert F(x,u,\alpha)-F(y,u,\beta)\vert \leq C_1(\vert x-y\vert +\vert \alpha-\beta\vert),
\end{equation}
for every $(x,y,\alpha,\beta,u)\in \R^d\times \R^d\times \R^k\times \R^k\times V$, 
%\item $G$ : $\R^d\rightarrow \R^{p\times k} $  is Lipschitz (and we denote by $C_2$ its Lipschitz constant)
\item {\bf{(H2)}} $A \in L^2((0,+\infty), M_{k\times d} )\cap L^1((0,+\infty), M_{k\times d})$ ($M_{k\times d}$ standing for the space of real  matrices with $k$ rows and $d$ columns).
\end{itemize}

In the sequel, we shall sometimes use a stronger assumption than  {\bf{(H2)}}. Namely:  {\bf{(H2')}} $A \in H^1((0,+\infty), M_{k\times d} )\cap L^1((0,+\infty), M_{k\times d})$.

\smallskip

Before studying the optimal control of equations with memory of type \pref{state}, let us establish the existence, uniqueness and continuous dependence with respect to initial conditions for  the Cauchy problem \pref{state}-\pref{bc}. The results of this section (Propositions \ref{cauchylip} and \ref{gron}) are fairly standard but we give proofs for the sake of completeness and to keep the present paper self-contained. Conditions  {\bf{(H1)}} and  {\bf{(H2)}} of course ensure existence and uniqueness of a solution to the Cauchy problem \pref{state}-\pref{bc}:

\begin{prop}\label{cauchylip}
Assume that  {\bf{(H1)}} and  {\bf{(H2)}} hold. For every $(x,z,u)\in \R^d\times L^2\times V$, the Cauchy problem \pref{state}-\pref{bc} admits a unique solution.
\end{prop}

\begin{proof}
For $\theta>0$, define
\[E_{\theta}:=\{y\in C^0(\R_+,\R^d),\; \sup_{t\geq 0} e^{-\theta t} \vert y(t)\vert <+\infty\}\]
and equip $E_\theta$ with the norm:
\[\Vert y\Vert_\theta :=  \sup_{t\geq 0} e^{-\theta t} \vert y(t)\vert.\]
Of course, $(E_\theta,\Vert .\Vert_\theta)$ is a Banach space. For $y\in E_{\theta}$, let us define:
\[Ty(t):=x+\int_0^t F(y(s),u(s),G_y(s))ds,\; \forall t\geq 0\]
where 
\[G_y(s):=\int_0^s A(\tau) y(s-\tau)d\tau+\int_s^{+\infty} A(\tau) z(\tau-s)d\tau.\] 
Until the end of the proof, $C$ will denote a positive constant (only depending on $F$ and $A$) which may vary from one line to another. Let $y\in E_\theta$, with our assumptions on $F$, we first get:
\begin{equation}\label{ineg0}
\vert Ty(t)\vert \leq \vert x \vert + C\left( t+\frac{e^{\theta t}}{\theta}\Vert y\Vert_{\theta}+\int_0^t\vert G_y(s)\vert ds  \right).
\end{equation}
Now, we also have
\[\begin{split}
\vert G_y(s)\vert &\leq \left( \int_0^s    \vert A(\tau)\vert  \vert y(s-\tau)\vert d\tau+ \Vert A \Vert_{L^2} \Vert z \Vert_{L^2}  \right)\\
&\leq   \Vert  A \Vert_{L^2}   \left(\Vert y\Vert_\theta  \left(\int_0^s  e^{2\theta(s- \tau) }d\tau  \right)^{1/2} + \Vert z \Vert_{L^2}  \right)\\
&\leq C\left(1+ \frac{e^{\theta s}\Vert y\Vert_{\theta}} {\sqrt{2\theta}}  \right).
\end{split}\]
Together with \pref{ineg0}, we then have
\[\vert Ty(t)\vert e^{-\theta t} \leq \vert x \vert e^{-\theta t} +C\left( te^{-\theta t}+ \Vert y\Vert_{\theta}\left(\frac{1}{\theta}+\frac{1}{\sqrt{2} \theta^{3/2}} \right)  \right)\]
which proves that $T(E_\theta)\subset E_\theta$.  For $y_1$ and $y_2$ in $E_\theta$ and $t\geq 0$, on the one hand, we have:
\[\vert Ty_1(t)-Ty_2(t)\vert \leq C\left( \frac{e^{\theta t}}{\theta}\Vert y_1-y_2\Vert_{\theta}+\int_0^t\vert G_{y_1}(s)-G_{y_2}(s)\vert ds  \right)\]
on the other hand:
\[\begin{split}
\vert G_{y_1}(s)-G_{y_2}(s)\vert 
 \leq  \frac{Ce^{\theta s}}{\sqrt{2\theta}} \Vert y_1-y_2 \Vert_\theta
\end{split}\]
so that:
\[\Vert Ty_1-Ty_2\Vert_{\theta}\leq C \Vert y_1-y_2 \Vert_\theta\left( \frac{1}{\theta}+\frac{1}{\sqrt{2} \theta^{3/2}}\right).\]
For $\theta$ large enough ($\theta\geq 2C+1$, say), $T$ is a contraction of $E_\theta$ hence admits a unique fixed-point. This clearly proves the desired result.
\end{proof}

From now on, for every $(x,z,u)\in \R^d\times L^2\times V$, we denote by $y_{x,z,u}$ the solution of  the Cauchy problem \pref{state}-\pref{bc}. The continuous dependence with respect to $(x,z)$ of trajectories of \pref{state}-\pref{bc} is given by:

\begin{prop}\label{gron}
Assume that  {\bf{(H1)}} and  {\bf{(H2)}} hold. Let $u\in V$, $(x_0,z_0)$ and $(x,z)$ be in $\R^d\times L^2$ and define $y_0:=y_{x_0,z_0,u}$, $y:=y_{x,z,u}$, then we have
\[\vert y(t)-y_0(t)\vert \leq Ce^{\theta t} (\vert x-x_0\vert + \Vert z-z_0 \Vert_{L^2}), \; \forall t\geq 0\]
for some constants $C$ and $\theta$ depending only on $F$ and $A$.

\end{prop}

\begin{proof}
In this proof, $C$ will denote a positive constant that only depends on $F$ and $A$ but which may vary from one line to another. Defining for $s\geq 0$
\[\begin{split}
\beta(s)&:=\int_0^s A(s-\tau)y(\tau)d\tau+\int_0^{+\infty} A(s+\tau)z(\tau) d\tau,\\
\beta_0(s)&:=\int_0^s A(s-\tau)y_0(\tau)d\tau+\int_0^{+\infty} A(s+\tau)z_0(\tau) d\tau,\\
\gamma(s)&:=\sup_{[0,s]} \vert y-y_0\vert,\; \Gamma(s):=\int_0^s \gamma,
\end{split}\]
we first have:
\[\vert y(s)-y_0(s)\vert \leq \vert x-x_0\vert + C\left(\int_0^s (\vert y-y_0\vert + \vert \beta-\beta_0\vert )\right)\]
Since we also have
\[\vert \beta(\tau)-\beta_0(\tau)\vert \leq C \left( \gamma(\tau)\Vert A \Vert_{L^1}+\Vert A \Vert_{L^2} \Vert z-z_0\Vert_{L^2}\right)\]
for $t\geq 0$, we then get:
\[\Gamma'(t)\leq \vert x-x_0\vert+  C \left(\Gamma(t)+ \Vert z-z_0\Vert_{L^2} t \right)\]
which, together with Gronwall's Lemma gives  the desired result. 
\end{proof}

\begin{Remark}\label{rem1}
Let us remark that when  one further assumes that {\bf{(H'2)}} holds (i.e. $A$ is further assumed to be $H^1$), then the estimate of proposition \ref{gron} also holds true when one replaces $ \Vert z-z_0 \Vert_{L^2}$ by  $\Vert z-z_0 \Vert_{(H^1)^{\prime}}$.
\end{Remark}
\subsection{The optimal control problem}

For $(x,z)\in \R^d\times L^2$, we consider the optimal control problem
\begin{equation}\label{defofv}
v(x,z):=\inf_{u\in V} \int_0^{+\infty} e^{-\lambda s} L(y_{x,z,u}(s),u(s))ds.
\end{equation}
where  {\bf{(H3)}}: $\lambda>0$ and $L$ : $\R^d\times K\rightarrow \R$ is assumed to be bounded, continuous and to satisfy
\begin{equation}\label{lipL}
\vert L(x,u)-L(y,u) \vert \leq C_2 \vert x-y\vert, \; \forall (x,y,u)\in \R^d\times \R^d\times K
\end{equation}
for some $C_2\geq 0$. Throughout the paper, we will assume that {\bf{(H1)}}, {\bf{(H2)}}  and  {\bf{(H3)}} hold.

\subsection{Continuity properties of the value function}

As a consequence of proposition \ref{gron}, we deduce that $v$ is bounded and uniformly continuous on $\R^d\times L^2$, which we denote $v\in \rm{BUC}(\R^d\times L^2,\R)$. More precisely, adapting classical arguments (see e.g. G. Barles \cite{barles}) to our context, we have:

\begin{prop}
Assume that  {\bf{(H1)}}, {\bf{(H2)}}  and  {\bf{(H3)}} hold, then  $v\in \rm{BUC}(\R^d\times L^2,\R)$ and more precisely, defining $\theta$ as in proposition \ref{gron}, one has:

\begin{enumerate}
\item $v$ is Lipschitz continuous on $\R^d\times L^2$ if $\lambda> \theta$,

\item $v\in C^{0,\alpha}(\R^d\times L^2,\R)$ for every $\alpha\in (0,1)$ if $\lambda= \theta$,

\item $v\in C^{0,\lambda/\theta}(\R^d\times L^2,\R)$ if $\lambda< \theta$.

\end{enumerate}

\end{prop}

\begin{proof}
Let us define
\[\delta:=\vert x-x_0\vert +\Vert z-z_0\Vert_{L^2}.\]
Let $\eps>0$ and $u_\eps$ be such that
\[\int_0^{+\infty} e^{-\lambda s} L(y_{x,z,u_\eps}(s),u_\eps(s))ds\leq v(x,z)+\eps\]
setting $y^{\eps}:=y_{x,z,u_\eps}$, $y_0^\eps:=y_{x_0,z_0,u_{\eps}}$ we then have:
\[v(x_0,z_0)-v(x,z)\leq\int_ 0^{+\infty} e^{-\lambda s}\left( L(y_0^\eps(s),u_\eps(s))-L(y^\eps(s),u_\eps(s))\right)ds+\eps.\]
Using proposition \ref{gron} and our assumptions on $L$, we then get, for some $C\geq 0$ and all $T\geq 0$:
\begin{equation}\label{ineg55}
v(x_0,z_0)-v(x,z)\leq C\left( \int_0^T \delta e^{(\theta-\lambda)s}ds +e^{-\lambda T}\right).\end{equation}
If $\lambda>\theta$, we then have:
\[v(x_0,z_0)-v(x,z)\leq \frac{C \delta }{\lambda -\theta}\]
which proves the first claim.

\smallskip

If $\lambda<\theta$ and if $\delta<1$ (which may be assumed to prove that $v$ is H\"older) taking $e^{-\lambda T}:=\delta^{\lambda/\theta}$ in \pref{ineg55} then yields
\[v(x_0,z_0)-v(x,z)\leq C\left(1+\frac{1}{\theta-\lambda} \right)  \delta^{\lambda/ \theta}\]
which proves the second claim.

\smallskip

Finally, if $\lambda=\theta$ (and again assuming $\delta<1$), taking $e^{-\lambda T}=\delta$ in  \pref{ineg55} yields:
\[v(x_0,z_0)-v(x,z)\leq  C  \left( \frac{-\delta \log(\delta)}{\lambda}+ \delta\right)\]
 which proves the last claim.

\end{proof}

\begin{Remark}\label{rem2}
Again, if $A$ is further assumed to be $H^1$ (i.e. when {\bf{(H'2)}} holds) then the uniform continuity of $v$ also holds true for the norm $ (x,z)\mapsto \vert x\vert + \Vert z \Vert_{(H^1)^{\prime}}$ i.e. when in the previous proof $\delta$ is replaced by $\delta:= \vert x-x_0\vert +\Vert z-z_0\Vert_{(H^1)^{\prime}}$. This fact will be useful later on when proving the comparison result. 
\end{Remark}

In the sequel, we shall denote by $C^0(\R^d\times \rm{L^2_w},\R)$ the class of real-valued functions defined on $\R^d\times L^2$ which are sequentially continuous for the weak topology of $\R^d\times L^2$, we then have the following:

\begin{prop}
Assume that  {\bf{(H1)}}, {\bf{(H2)}}  and  {\bf{(H3)}} hold, then $v\in C^0(\R^d\times \rm{L^2_w},\R)$.
\end{prop}

\begin{proof}
Let $(\alpha_n)_n :=(z_n,x_n)_n$ be a weakly convergent  sequence in $\R^d\times L^2$ and let us denote  by $\alpha:=(x,z)\in \R^d\times L^2$ its weak limit. Let $u\in V$ be some admissible control and simply denote $y_n:=y_{\alpha_n, u}$ and $y:=y_{\alpha, u}$ the trajectories of \pref{state} associated respectively to the initial conditions $\alpha_n$ and $\alpha$.  If we prove that $y_n$ converges uniformly on compact subsets to $y$ as $n$ tends to $+\infty$ then the desired result will easily follow from our assumptions on $L$. Let us define
\[\delta_n(t):=\int_0^{+\infty} A(t+s) z_n(s) ds, \; \delta(t):= \int_0^{+\infty} A(t+s) z(s) ds.\]
Since $A\in L^2$ by {\bf{(H2)}}, we have:
\[\vert \delta_n (t) \vert \leq \Vert A\Vert_{L^2}  \Vert z_n \Vert_{L^2}\leq C.\]
Thanks to {\bf{(H2)}} again, $\delta_n$ converges pointwise to $\delta$. Rewriting the state equation as:
\[\begin{split}
 \dot{y}(t) &= F\left(y(t),u(t), \delta(t)+\int_0^t A(t-s)y(s)ds\right),\\
 \dot{y}_n(t) & = F\left(y_n(t),u(t), \delta_n(t)+\int_0^t A(t-s)y_n(s)ds\right)\end{split}\]
 we get:
 \begin{equation}\label{inegtraj1}
 \vert \dot{y}_n-\dot{y}\vert(t) \leq C \left (\vert y_n-y\vert(t) +\vert \delta_n-\delta\vert(t) + \int_0^t \vert   A(t-s) (y_n-y)(s)\vert ds   \right)
 \end{equation}
 (where again in this proof $C$ denotes a nonnegative constant depending only on $F$ and $A$ but possibly changing from one line to another). Defining 
 \[ \gamma_n(t):= \sup_{[0,t]} \vert y_n-y \vert,\; \Gamma_n(t):=\int_0^t \gamma_n,\]
inequality \pref{inegtraj1} yields for all $s\in [0,t]$:
 \begin{equation*}
 \vert \dot{y}_n-\dot{y}\vert(s) \leq C \left(\vert y_n-y\vert(s) +\vert \delta_n-\delta\vert(s) + \Vert A\Vert_{L^1(0,t)} \gamma_n(s)  \right).
 \end{equation*}
Integrating the previous yields:
\[\vert y_n- y\vert(s) \leq \vert x_n -x \vert +C \left ( \Gamma_n(t)+ \int_0^t \vert \delta_n -\delta \vert  \right),\; \forall s\in [0,t].\]
Hence 
\begin{equation}\label{inegtraj2}
\gamma_n(t)= \dot{\Gamma}_n(t)\leq \vert x_n -x \vert + C\left( \Gamma_n(t)+ \int_0^t \vert \delta_n -\delta \vert  \right).
 \end{equation}
On the one hand, Dominated convergence implies that 
\[ \lim_n \int_0^{+\infty} \vert \delta_n -\delta \vert =0,\]
on the other hand, \pref{inegtraj2} and Gronwall's Lemma  imply that $\Gamma_n(t)$ tends to $0$ as  $n$ tends to $+\infty$. With  \pref{inegtraj2}, this proves that $\gamma_n(t)$ tends to $0$ as  $n$ tends to $+\infty$ and the desired result follows. Let us remark that, with \pref{inegtraj1}, this of course also implies that $(y_n)_n$ converges to $y$ in $W^{1,\infty}_{\rm{loc}}(\R_+,\R^d)$.

\end{proof}

\section{Dynamic programming principle}\label{sec2}

Our aim now is to prove that the value function $v$ obeys the following dynamic programming principle:

\begin{prop}\label{progd}
Let $(x,z)\in \R^d\times L^2$ and $t\geq 0$, we then have:
\begin{equation}\label{ppd}
v(x,z)=\inf_{u\in V} \left\{  \int_0^t e^{-\lambda s}  L(y_{x,z,u}(s),u(s))ds + e^{-\lambda t} v(y_{x,z,u}(t), y_{x,z,u}(t-.)) \right\}
\end{equation}
$(y_{x,z,u}(t-.)(s):=y_{x,z,u}(t-s)$ for all $s>0)$.
\end{prop}

\begin{proof}
Let $\eps>0$ and $u_{\eps}\in V$ be such that
\[ \int_0^{+\infty} e^{-\lambda s}  L(y_{x,z,u_\eps}(s),u_\eps(s))ds \leq v(x,z)+\eps\]
we then have
\[\begin{split} v(x,z)+\eps \geq \int_0^t e^{-\lambda s}  L(y_{x,z,u_\eps}(s),u_\eps(s))ds\\
 +e^{-\lambda t} \int_0^{+\infty}\
  e^{-\lambda \tau} L(y_{x,z,u_{\eps}}(t+\tau), u_{\eps}(t+\tau))d\tau.\end{split}\]
Using the fact that  $y_{x,z,u_{\eps}}(t+.)$ is the trajectory associated to the initial conditions $(y_{x,z,u_{\eps}}(t)$,  $y_{x,z,u_{\eps}}(t-.))$ and the control $u_{\eps}(t+.)$, we deduce:
\[\begin{split}
 v(x,z)+\eps &\geq    \int_0^t e^{-\lambda s}  L(y_{x,z,u_\eps}(s),u_\eps(s))ds+e^{-\lambda t} v(y_{x,z,u_\eps}(t), y_{x,z,u_\eps}(t-.))\\
 &\geq \inf_{u\in V} \left\{  \int_0^t e^{-\lambda s}  L(y_{x,z,u}(s),u(s))ds + e^{-\lambda t} v(y_{x,z,u}(t), y_{x,z,u}(t-.)) \right\}.
 \end{split}\]

To prove the converse inequality, let $u\in V$, 
\[x_t:= y_{x,z,u}(t),\;  z_t:=y_{x,z,u}(t-.),\]
$\eps>0$ and $\omega_{\eps}\in V$ be such that:
\[\int_0^{+\infty} e^{-\lambda s}  L(y_{x_t,z_t,\omega_\eps}(s),\omega_\eps(s))ds \leq v(x_t,z_t)+\eps.\]
defining 
\[u_{\eps}(s):= \left\{
\begin{array}{lll}
u(s) &  \mbox{ if } s\in[0,t]\\
\omega_{\eps}(s-t)  & \mbox{ if } s>t
\end{array}\right.\]
we have
\[y_{x,z,u_\eps}(s):= \left\{
\begin{array}{lll}
y_{x,z,u}(s) &  \mbox{ if } s\in[0,t]\\
y_{x_t,z_t,\omega_{\eps}}(s-t)  & \mbox{ if } s>t
\end{array}\right.\]
hence
\[\begin{split}
v(x,z) &\leq  \int_0^t e^{-\lambda s}  L(y_{x,z,u}(s),u(s))ds +\int_t^{+\infty} e^{-\lambda s}  L(y_{x_t,z_t,\omega_{\eps}}(s-t),\omega_\eps(s-t))ds\\
&\leq \int_0^t e^{-\lambda s}  L(y_{x,z,u}(s),u(s))ds+e^{-\lambda t} (v(x_t,z_t)+\eps)
\end{split}\]
since $u$ and $\eps>0$ are arbitrary, we get the desired result.

\end{proof}

 \section{Viscosity solutions and Hamilton-Jacobi-\\
 Bellman equations}\label{sec3}

\subsection{Preliminaries}

Let us define
\begin{equation}
E_0:=\{(z(0),z) \mbox{, } z\in H^1((0,+\infty),\R^d)\}
\end{equation}
and remark that $E_0$ is a dense subspace of our initial state space $E:=\R^d\times L^2$. With the uniform continuity of $v$ on $\R^d\times L^2$, this implies that $v$ is fully determined by its restriction to $E_0$. In fact, we will derive from the dynamic programming principle a PDE satisfied by $v$ (a priori only on $E_0$) and by a comparison result we will see that in fact this equation (satisfied in some appropriate  viscosity sense) fully characterizes the value function $v$. Let us define
\begin{equation}\label{defF}
\F(x,z,u):=F\left(x,u,\int_0^{+\infty} A(s) z(s) ds\right),\; \forall (x,z,u)\in\R^d\times L^2\times K.
\end{equation}

Before going further, we need the following classical lemma:

\begin{lem}\label{transl}
Let $\delta>0$ and $z\in L^2((-\delta,+\infty),\R^d)$ for $t\in [0, \delta]$, define $z_t(s):=z(s-t)$ for $s\geq 0$, then $z_t$ converges to $z$ in $L^2(\R_+,\R^d)$ as $t$ goes to $0^+$.
\end{lem}

\begin{proof}
First, let us remark that $\Vert z_t\Vert_{L^2(\R_+)}$ converges to $\Vert z \Vert_{L^2(\R_+)}$ hence it is enough to prove that $z_t$ converges weakly to $z$ in $L^2(\R_+,\R^d)$. To see this, let $f\in C_c^0(\R_+,\R^d)$, and remark that
\[\int_{\R_+} (z_t-z)\cdot f=\int_{-t}^0 z(s)\cdot f(s+t)ds+\int_{\R_+} z(s)\cdot (f(s+t)-f(s))ds\]
Since both terms in the right-hand side tend to $0$ as $t$ tends to $0$, this proves  that $z$ is the only weak limit point of the bounded sequence $z_t$ hence that $z_t$ converges weakly to $z$ in $L^2(\R_+,\R^d)$. 

\end{proof}

We will also need the following

\begin{lem}\label{derivtp}
Let $(x,z)\in E_0$, $u\in V$ and for $t>0$ define:
\[x_t:=y_{x,z,u}(t),\; z_t(\tau):=y_{x,z,u}(t-\tau), \; \forall \tau>0,\]
then $t\mapsto (x_t,z_t)$ is locally Lipschitz in $t$ (uniformly in the control $u\in V$) for the $\R^d\times L^2$ norm. Moreover, for all $t\geq 0$
\begin{equation}\label{derivpast}
\lim_{s\rightarrow 0^+} \frac{z_{t+s}-z_t}{s}=-\dot{z}_t \mbox{ in } L^2,
\end{equation}
and, for every $t$ which is a Lebesgue point of $t\mapsto \F(x_t,  z_t,u(t))$ 
\begin{equation}\label{deriveta}
\lim_{s\rightarrow 0^+} \frac{x_{t+s}-x_t}{s}=\F(x_t,  z_t,u(t)).
\end{equation}
\end{lem}

\begin{proof}
The  lipschitzianity of $t\mapsto x_t$ and the proof of \pref{deriveta} are straightforward. To shorten notation, we define $y:=y_{x,z,u}$ and remark that for every $t\geq 0$, $y\in H^1((-\infty,t),\R^d)$ and that $y\in W^{1,\infty}_{\rm{loc}}((0,+\infty),\R^d)$ so that \pref{derivpast} will imply the  local lipschitzianity of $t\mapsto z_t$. To prove  \pref{derivpast}, let us introduce for $s>0$ and $\tau\geq 0$:
\[\Delta_s(\tau):=\frac{z_{t+s}(\tau)-z_t(\tau)}{s}+\dot{z}_t(\tau)\]
which can be rewritten as:
\[\begin{split}
\Delta_s(\tau)&=\frac{y(t+s-\tau)-y(t-\tau)}{s}-\dot{y}(t-\tau)\\
&=\frac{1}{s}  \int_0^s \left(\dot{y}(t-\tau+\alpha)-\dot{y}(t-\tau)  \right)  d\alpha\end{split}\]
Jensen's inequality first yields:
\[\Delta_s(\tau)^2\leq \frac{1}{s}  \int_0^s \left(\dot{y}(t-\tau+\alpha)-\dot{y}(t-\tau)  \right)^2  d\alpha\]
using Fubini's theorem, we then  get:
\[\Vert \Delta_s \Vert^2_{L^2}\leq \frac{1}{s}  \int_0^s \left( \int_{\R_+} \left (\dot{y}(t-\tau+\alpha)-\dot{y}(t-\tau)  \right)^2 d\tau\right) d\alpha.\]
By the same arguments as in lemma \ref{transl} and since $y\in H^1((-\infty,t),\R^d)$, for every $t\geq 0$, we deduce that for every $\eps>0$, there exists $s_\eps>0$ such that for all $\alpha\in[0,s_{\eps}]$, one has:
\[\left( \int_{\R_+} \left (\dot{y}(t-\tau+\alpha)-\dot{y}(t-\tau)  \right)^2 d\tau\right)\leq \eps\]
which proves that $\Delta_s$ converges to $0$ in $L^2$ as $s$ goes to $0$.  
\end{proof}

\begin{Remark}\label{rem3}
It follows from lemma \ref{derivtp} that if $\psi\in C^1(\R\times\R^d\times L^2,\R)$ then the function $t\mapsto \psi(t,x_t,z_t)$ is locally Lipschitz in $t$ (in fact, uniformly in the control $u\in V$) and for all $t\geq 0$, one has:
\[\begin{split}
&\psi(t,x_t,z_t)=\psi(0,x_0,z_0)\\
&+\int_0^t \left( \partial_t \psi(s,x_s,z_s)+\nabla_x \psi(s,x_s,z_s)\cdot \F(x_s,z_s,u(s))-\<D_z\psi(s,x_s,z_s), \dot{z}_s>  \right)ds. \end{split}\]
\end{Remark}

\subsection{Formal derivation of the equation}

To formally derive the  Hamilton-Jacobi-Bellman equation of our problem, let us assume for a moment that $v$ is of class $C^1$ on $\R^d\times L^2$. Let $(x,z)\in E_0$, $u\in K$ be some admissible constant control, and set 
\[(x_t,z_t):=(y_{x,z,u}(t), y_{x,z,u}(t-.)),\]
 by the dynamic programming principle, we first have:
\[v(x,z)\leq \int_0^t e^{-\lambda s} L(x_s, u)ds +e^{-\lambda t} v(x_t,z_t)\]
so that
\[L(x,u)+ \lim_{t\rightarrow 0^+} \frac{1}{t} \left(  e^{-\lambda t} v(x_t,z_t)-v(x,z) \right)\geq 0\]
together with lemma \ref{derivtp}, this reads as
\[L(x,u)+\nabla_x v(x,z)\cdot \F(x,z,u) -\lambda v(x,z)-\<D_z v(x,z), \dot{z}>\geq 0\]
and since $u$ is arbitrary, this yields
\[\lambda v(x,z)+\<D_z v(x,z),\dot{z}> + \sup_{u\in K} \{-L(x,u)- \nabla_x v(x,z)\cdot \F(x,z,u)\}\leq 0.\] 

We then define the Hamiltonian:
\begin{equation}\label{defH}
H(x,z,p):=\sup_{u\in K} \{-L(x,u)-p\cdot \F(x,z,u)\},\; \forall (x,z,p)\in\R^d\times L^2\times \R^d.
\end{equation}

Let $u\in V$, and simply denote $(x_{t,u},z_{t,u}):=(y_{x,z,u}(t), y_{x,z,u}(t-.))$. Using  remark \ref{rem3} following lemma \ref{derivtp}, with $\psi(t,x,z):=e^{-\lambda t} v(x,z)$, we have:
\[\begin{split}
&e^{-\lambda t}v(x_{t,u},z_{t,u})=v(x,z)-\int_0^t e^{-\lambda s} \lambda v(x_{s,u}, z_{s,u})ds\\
&+\int_0^t e^{-\lambda s} \left(\nabla_x v(x_{s,u},z_{s,u})\cdot \F(x_{s,u},z_{s,u},u(s))
-\<D_z v(x_{s,u},z_{s,u}),\dot{z}_{s,u}> \right) ds 
\end{split}\]
so that the dynamic programming principle yields:
\[\begin{split}
&0=  \inf_{u\in V} \{ \int_0^t  e^{-\lambda s}( L(x_{s,u},u(s))- \lambda v(x_{s,u}, z_{s,u})
+ \nabla_x v(x_{s,u},z_{s,u})\cdot \F(x_{s,u},z_{s,u},u(s)) \\
& -\<D_z v(x_{s,u},z_{s,u}),\dot{z}_{s,u}> ) ds  \} \\
&\geq \inf_{u\in V} \{ \int_0^t e^{-\lambda s} ( -H(x_{s,u},z_{s,u}, \nabla_x v(x_{s,u},z_{s,u}))- \lambda v(x_{s,u}, z_{s,u})\\
&-\<D_z v(x_{s,u},z_{s,u}),\dot{z}_{s,u}> ) ds \}.
\end{split}\]
It is natural to expect the integrand above to converge  as $t\rightarrow 0^+$, uniformly in $u$ to  
\[-H(x,z, \nabla_x v(x,z))- \lambda v(x, z) -\<D_z v(x,z),\dot{z}>\]
so that:
\[\lambda v(x,z)+\<D_z v(x,z),\dot{z}> + H(x,z, \nabla_x v(x,z))\geq 0.\] 

\smallskip

Thus,  at least formally, the Hamilton-Jacobi-Bellman equation satisfied by the value function $v$ can be written as:
\begin{equation}\label{hjb}
\lambda v(x,z)+H(x,z,\nabla_x v(x,z))+\<D_z v(x,z),\dot{z}>=0 \mbox{ on } E_0.
\end{equation}
with $H$ defined by \pref{defH}.

\smallskip

The next lemma whose easy proof is left to the reader gives the regularity properties of $H$ 
\begin{lem}\label{modulh}
Let $H$ be the Hamiltonian defined by \pref{defH}. Assume that  {\bf{(H1)}}, {\bf{(H2)}}  and  {\bf{(H3)}} hold, then  there exists a nonnegative constant $C$ such that:
 \begin{equation}\label{ucdeH10}
 \vert H(x,z,p)-H(y,w,p)\vert \leq C\left(\vert x-y\vert + \Vert z-w\Vert_{L^2} \right) (1+\vert p\vert),  \end{equation}
 and 
  \begin{equation}\label{ucdeH20}
 \vert H(x,z,p)-H(x,z,q)\vert \leq C \vert p-q \vert  (1+ \vert x\vert + \Vert z \Vert_{L^2}  ), 
 \end{equation}
 for every $(x,z, y, w ,p, q)\in (\R^d\times L^2)^2\times \R^d\times \R^d$. If, in addition {\bf{(H'2)}} is satisfied, then \pref{ucdeH10} can be improved by:
 \begin{equation}\label{ucdeH100}
 \vert H(x,z,p)-H(y,w,p)\vert \leq C\left(\vert x-y\vert + \Vert z-w\Vert_{(H^1)^{\prime}} \right) (1+\vert p\vert),
 \end{equation}
for every  $(x,z,y,w,p)\in (\R^d\times L^2)^2 \times \R^d$.
\end{lem}

\subsection{Definition of viscosity solutions}

The formal manipulations above actually suggest that the natural definition of viscosity solutions in the present context  should read as:
\begin{defi}\label{defvisc}
Let $w\in \rm{BUC}(\R^d\times L^2,\R)\cap C^0(\R^d\times L^2_w,\R)$, then $w$ is said to be
\begin{enumerate}

\item a viscosity subsolution of \pref{hjb} on $\R^d\times L^2$  if for every $(x_0,z_0)\in \R^d\times L^2$ and every $\phi\in C^1(\R^d\times L^2,\R)$ such that $w-\phi$ has a local maximum (in the sense of the strong topology of $\R^d\times L^2$)
 at $(x_0,z_0)$, one has:
 \[\lambda w(x_0,z_0)+H(x_0,z_0,\nabla_x \phi(x_0,z_0))+\mathrm{liminf}_{(x,z)\in E_0 \rightarrow (x_0,z_0) } \<D_z \phi(x,z),\dot{z}>\leq 0,\]
 where the convergence $(x,z)\in E_0 \rightarrow (x_0,z_0)$ has to be understood in the strong $\R^d\times L^2$ sense,
 \item a viscosity supersolution of \pref{hjb} on  $\R^d\times L^2$  if for every $(x_0,z_0)\in \R^d\times L^2$ and every $\phi\in C^1(\R^d\times L^2,\R)$ such that $w-\phi$ has a local minimum (in the sense of  the strong topology of $\R^d\times L^2$)
 at $(x_0,z_0)$, one has:
 \[\lambda w(x_0,z_0)+H(x_0,z_0,\nabla_x \phi(x_0,z_0))+\mathrm{limsup}_{(x,z)\in E_0 \rightarrow (x_0,z_0)} \<D_z \phi(x,z),\dot{z}>\geq 0,\]
 where the convergence $(x,z)\in E_0 \rightarrow (x_0,z_0)$ has to be understood in the strong $\R^d\times L^2$ sense,

\item a viscosity solution of \pref{hjb} on $\R^d\times L^2$  if it is both a  viscosity subsolution of \pref{hjb} and a viscosity supersolution of \pref{hjb} on $\R^d\times L^2$.

\end{enumerate}

\end{defi}

\subsection{The value function is a viscosity solution}

\begin{prop}\label{vvisc}
Assume that  {\bf{(H1)}}, {\bf{(H2)}}  and  {\bf{(H3)}} hold. The value function $v$ defined by \pref{defofv} is a viscosity solution of \pref{hjb} on $\R^d\times L^2$.
\end{prop}

\begin{proof}
{\bf{Step 1}}: $v$ is a viscosity subsolution.

\medskip

Let $\alpha_0:=(x_0,z_0)\in \R^d\times L^2$ and $\phi\in C^1(\R^d\times L^2,\R)$ such that $v(x_0,z_0)=\phi(x_0,z_0)$ and $\phi\geq v$ on the ball $B_r:=B((x_0,z_0), r)$ of $\R^d\times L^2$. 
Let $u\in K$ be some constant control. For all $\eps>0$, there is some $\alpha_{\eps}:=(x_{\eps},z_{\eps})\in E_0$ such that
\[\phi(\alpha_{\eps})-\eps^2\leq  v(\alpha_{\eps})\leq \phi(\alpha_\eps),\; \lim_{\eps} \alpha_\eps= \alpha_0 \mbox{ in } \R^d \times L^2\]
and such that $\alpha_{\eps, s}= (x_{\eps,s}, z_{\eps, s})= (y_{\alpha_{\eps}, u }(s),  y_{\alpha_{\eps}, u }(s-.))$ belongs to $B_r$ for all $s\in [0,\eps]$.  Note that by construction, $\alpha_{\eps,s}$ belongs to $E_0$ for every $s\in [0,\eps]$. The dynamic programming principle first yields
\begin{equation}\label{inegss1}
\phi(\alpha_{\eps})-\eps^2 \leq \int_0^\eps   e^{-\lambda s} L(x_{\eps,s},u)ds+e^{-\lambda \eps} \phi(\alpha_{\eps,\eps}).\end{equation}
Thanks to the smoothness of $\phi$,  lemma \ref{derivtp} and  remark \ref{rem3}, we can write:
\[\begin{split}
e^{-\lambda \eps} \phi(\alpha_{\eps,\eps})=\phi(\alpha_{\eps})
+\eps(\nabla_x \phi(\alpha_{\eps})\cdot \F(\alpha_{\eps},u)-\lambda \phi(\alpha_\eps))\\
-\int_0^{\eps} e^{-\lambda s} \<D_z \phi(\alpha_{\eps,s}),\dot{z}_{\eps,s}>  ds   +o(\eps)\end{split}\]

Using \pref{inegss1}, dividing by $\eps$ and taking the liminf as $\eps \to 0^+$, we then get:
\[\begin{split}
0  &\geq  \mathrm{liminf}_{\eps \to 0^+} \left(  - \frac{1}{\eps}  \int_0^{\eps} e^{-\lambda s} L(x_{\eps,s},u)ds-\nabla_x \phi(\alpha_{\eps})\cdot \F(\alpha_{\eps},u)+\lambda \phi(\alpha_\eps) \right)\\
&+   \mathrm{liminf}_{\eps \to 0^+} \frac{1}{\eps}  \int_0^{\eps} e^{-\lambda s}  \<D_z \phi(\alpha_{\eps,s}),\dot{z}_{\eps,s}> ds 
\end{split}.\] 
Since we have:
\[  \mathrm{liminf}_{\eps \to 0^+} \frac{1}{\eps}  \int_0^{\eps} e^{-\lambda s}  \<D_z \phi(\alpha_{\eps,s}),\dot{z}_{\eps,s}> ds 
\geq \mathrm{liminf}_{(x,z)\in E_0\to \alpha_0} \<D_z \phi(x,z),\dot{z}>, \]
we then obtain:
\[0 \geq -L(x_0,u)-\nabla_x \phi(\alpha_{0})\cdot \F(\alpha_{0},u)+\lambda v(\alpha_0)+
\mathrm{liminf}_{(x,z)\in E_0\to \alpha_0} \<D_z \phi(x,z),\dot{z}>.\]
Since $u\in K$ is arbitrary in the previous inequality, taking the supremum with respect to $u$ and using the very definition of $H$ given in \pref{defH}, we thus deduce:
 \[\lambda v(\alpha_0)+H(\alpha_0,\nabla_x \phi(\alpha_0))
 +\mathrm{liminf}_{(x,z)\in E_0\rightarrow \alpha_0} \<D_z \phi(x,z),\dot{z}>\leq 0,\]
which proves that  $v$ is a viscosity subsolution of \pref{hjb} on $\R^d\times L^2$.

\smallskip

{\bf{Step 2}}: $v$ is a viscosity supersolution.

\medskip

Now let $\phi\in C^1(\R^d\times L^2,\R)$ be such that $v(\alpha_0)=\phi(\alpha_0)$ and $v \geq \phi$ on the ball $B_r:=B(\alpha_0, r)$ of $\R^d\times L^2$. For all $\eps>0$, there is some $\alpha_{\eps}:=(x_{\eps},z_{\eps})\in E_0$ such that
\[\phi(\alpha_{\eps})+\eps^2\geq  v(\alpha_{\eps})\geq \phi(\alpha_\eps)\]
and such that for every $u\in V$,  $\alpha_{\eps, s,u}= (x_{\eps,s,u}, z_{\eps, s,u})= (y_{\alpha_{\eps}, u }(s),  y_{\alpha_{\eps}, u }(s-.))$ belongs to $B_r$ for all $s\in [0,\eps]$. 
The dynamic programming principle then gives:
\begin{equation}\label{inegsursol}
\phi(\alpha_\eps)+\eps^2 \geq \inf_{u\in V} \left\{ \int_0^\eps e^{-\lambda s} L(x_{\eps, s,u},u(s))ds+e^{-\lambda \eps} \phi(\alpha_{\eps,\eps,u})\right\}\end{equation}
With lemma \ref{derivtp} and  remark \ref{rem3}, we can rewrite
\[\begin{split}
&e^{-\lambda \eps}\phi(\alpha_{\eps,\eps,u})=\phi(\alpha_\eps)-\int_0^\eps e^{-\lambda s} \lambda \phi(\alpha_{\eps,s,u})ds\\
&+\int_0^\eps e^{-\lambda s} \left(\nabla_x \phi(\alpha_{\eps,s,u})\cdot \F(\alpha_{\eps,s,u},u(s))
-\<D_z \phi(\alpha_{\eps, s,u}),\dot{z}_{\eps,s,u}> \right) ds 
\end{split}\]
Using \pref{inegsursol}, we then get:
\[\begin{split}
&\eps^2\geq  \inf_{u\in V} \{ \int_0^\eps   e^{-\lambda s}( L(x_{\eps,s,u},u(s))- \lambda \phi(\alpha_{\eps,s,u})
+ \nabla_x \phi(\alpha_{\eps,s,u})\cdot \F(\alpha_{\eps,s,u},u(s)) \\
& -\<D_z \phi(\alpha_{\eps,s,u}),\dot{z}_{\eps, s,u}> ) ds  \} \\
&\geq \inf_{u\in V} \{ \int_0^\eps e^{-\lambda s} ( -H(\alpha_{\eps,s,u}, \nabla_x \phi(\alpha_{\eps,s,u}))- \lambda \phi(\alpha_{\eps,s,u})-\<D_z \phi(\alpha_{\eps,s,u},\dot{z}_{\eps, s,u}> ) ds \}
\end{split}\]
From the continuity of  $\nabla_x \phi$, $D_z \phi$, and $H$, we deduce the (uniform in $u$) convergence as $\eps \rightarrow 0^+$ and $s\in [0,\eps]$, of $\nabla_x \phi(\alpha_{\eps,s,u})$,  $D_z \phi(\alpha_{\eps, s,u})$ and \\ 
$H(\alpha_{\eps,s,u}, \nabla_x \phi(\alpha_{\eps,s,u}))$ respectively to $\nabla_x \phi(\alpha_0)$, $D_z \phi(\alpha_0)$ and 
$H(\alpha_{0}, \nabla_x \phi(\alpha_{0}))$. Dividing by $-\eps$ the last inequality and taking the limsup as $\eps\to 0^+$, we thus get:
\[  \lambda v(\alpha_{0})+ H(\alpha_{0}, \nabla_x \phi(\alpha_{0}))
+\mathrm{limsup}_{(x,z)\in E_0\rightarrow \alpha_0}  \<D_z \phi(x,z),\dot{z}>\geq 0\]
which proves that  $v$ is a viscosity supersolution of \pref{hjb} on $\R^d\times L^2$.

\end{proof}

\section{Comparison and uniqueness}\label{sec4}

\subsection{Preliminaries}

Our aim now is to prove that $v$ is the unique viscosity subsolution of the Hamilton-Jacobi-Bellman equation:
\begin{equation*}
\lambda v(x,z)+H(x,z,\nabla_x v(x,z))+\<D_z v(x,z),\dot{z}>=0.
\end{equation*}
This will of course follow from a comparison result stating that if $v_1$ and $v_2$ (in a suitable class of continuous functions) are respectively a viscosity subsolution and a viscosity supersolution of the equation then $v_1\leq v_2$ on $E:=\R^d\times L^2$. As usual, the comparison result is proved, by introducing a doubling of variables and by considering  perturbed problems of the form:
\[ \sup \left\{ v_1(\alpha_1)-v_2(\alpha_2)- P_\theta(\alpha_1,\alpha_2)\; : \; (\alpha_1, \alpha_2) \in (\R^d\times L^2)^2Ê\right\}\]
where $P_\theta$ is some  perturbation function (depending on small parameters $\theta$). This perturbation  includes a penalization of the doubling of variables and coercive terms that ensure the existence of maxima say $\alpha_1^\theta$ and $\alpha_2^\theta$. Then one uses the fact that $v_1$ is a viscosity subsolution by taking $\phi:=P_\theta (.,\alpha_2^\theta)$  as test-function.

\smallskip

To overcome the difficulties due to the fact  that the term $\<D_z \phi(\alpha),\dot{z}>$  is only defined for $z\in H^1$ and that the equation is only justified when in addition $z(0)=x$, one has to be careful on the choice of the perturbation $P_\theta$. For general infinite-dimensional Hamilton-Jacobi equations with an unbounded linear term, these difficulties were solved in a general way by M. Crandall and P.-L. Lions in  \cite{cranlio4, cranlio5, cranlio6}. One of the key arguments of M. Crandall and P.-L. Lions in \cite{cranlio4} is to use a suitable norm to penalize the doubling of variables. In our context this, roughly speaking, amounts to use \emph{a kind of}  $(H^1)^{\prime}$ norm instead of the  $L^2$ norm in the doubling of variables.  As usual, the comparison proof will very much rely on the use of quadratic test-functions of the form:
\[ \phi (\alpha)= \Phi_1(x)+\Phi_2(\alpha)+\Phi_3(z):= a\vert x \vert^2 + b\<B(\alpha-\alpha_0),\alpha-\alpha_0>+c\Vert z \Vert^2\]
where $a$ and $b$ and $c$ are constants and $B$ is a bounded positive self-adjoint operator of $\R^d\times L^2$. Let us first note that in our case the term $\< \Phi'_3(z),\dot{z}>=2c\<z,\dot{z}>$ can be dealt easily since, if $z\in H^1$, one has:
\[\<z,\dot{z}>=-\frac{1}{2} \vert z(0)\vert^2\]
and since in the definition of viscosity solutions, we have imposed the convergence of the initial value $z(0)$, it is easy to figure out that this term won't be a big problem in the proof. The difficulty of dealing with the second term $\< D_z \Phi_2(\alpha),\dot{z}>$  can be solved by properly choosing $B$ as in M. Crandall and P.-L. Lions \cite{cranlio5} who emphasized the good properties $B$ should enjoy for the comparison proof to work.  We now proceed to the explicit construction of such a $B$ in our context.

\smallskip

Let us endow $E:=\R^d\times L^2$ with its standard Hilbertian structure, i.e. with the norm: 
\[\Vert \alpha \Vert^2:= \vert x\vert^2 +\Vert z \Vert_{L^2}^2, \; \forall \alpha=(x,z)\in E\]
and the corresponding inner product $\<.,.>$. Let $T$ be the linear unbounded operator on $E$ with domain $D(T)=\R^d\times H^1$ and defined by
\begin{equation}\label{defofT}
T(y,w):=(y-w(0),-\dot{w}),\; \forall (y,w)\in D(T).
\end{equation}
Its adjoint $T^*$ has domain $D(T^*)=E_0=\{(x,z)\in E \; : \; z\in H^1,\; z(0)=x\}$ and is given by
\begin{equation}
T^*(x,z):=(z(0),\dot{z})=(x,\dot{z}),  \; \forall (x,z)\in D(T^*)=E_0.
\end{equation}\label{defofTad}
The unbounded operator $I+T^*T$ therefore has domain 
\[D(T^*T)=\{(y,w)\in E\; : \; w \in H^2, \; y=w(0)-\dot{w}(0)\}\]
and is given by
\begin{equation}
(I+T^*T)(y,w):=(2y-w(0),-\ddot{w}+w),  \; \forall (y,w)\in D(T^*T).
\end{equation}
Now, let us set $B:=(I+T^*T)^{-1}$. For $\alpha=(x,z)\in E$, $(y,w):=B(\alpha)$ is defined as follows:  firtsly, $w\in H^2$ is the solution of 
\begin{equation}\label{defBw}
\left\{\begin{array}{cccc}
-\ddot{w} + w    & = &  z  & \mbox{ in } (0,+\infty)\\
-2 \dot{w}(0) + w(0)   & = & x, & 
\end{array}\right.
\end{equation}
secondly, $y$ is defined by 
\begin{equation}\label{defBy}
y=\frac{x+w(0)}{2}=w(0)-\dot{w}(0).
\end{equation}
In the sequel we shall also denote $B=(B_1,B_2)$ where $B_2(x,z)=w$ is defined by \pref{defBw} and $B_1(x,z)=y$ is given by \pref{defBw} and \pref{defBy}.
Setting:
\[\Vert \alpha \Vert_B^2:= \<B(\alpha),\alpha>,\; \forall \alpha=(x,z)\in E,\] 
and defining $w:=B_2(\alpha)$ by \pref{defBw}, an elementary computation shows that
\begin{equation}\label{balal}
\Vert \alpha \Vert_B^2= \frac{\vert x\vert ^2}{2}+\frac{\vert w(0)\vert^2}{2}+\Vert w \Vert_{H^1}^2\end{equation}
which in particular shows that there is some constant $C>0$, such that
\begin{equation}\label{controlhprime}
\Vert (x,z)\Vert_B^2\geq C\left (\vert x\vert^2+ \Vert z \Vert^2_{(H^1)^{\prime}}\right),\; \forall (x,z)\in E.
\end{equation}
Obviously, by construction $B$ is  a self-adjoint, nonnegative compact operator on $E$ and $TB$ is a bounded operator on $E$. For $\alpha=(x,z)\in E_0=D(T^*)$ and $w:=B_2(\alpha)$, some computations lead to:
\[\<TB(\alpha),\alpha>=\<B(\alpha),T^*(\alpha)>=\frac{3}{8}\vert x-w(0)\vert^2+\frac{ x \cdot w(0)}{2} \geq \frac{1}{8}\vert x-w(0)\vert^2,\]
and since $TB$ is continuous and $E_0$ is dense in $E$ this proves
\begin{equation}\label{tbpositif}
\<TB(\alpha),\alpha>\geq 0, \; \forall \alpha\in E.
\end{equation} 
Let us also remark that for $\alpha\in E$ and $\beta:=B(\alpha)$ one has:
\[\Vert \alpha\Vert_B^2=\<\alpha, \beta>=\<(I+T^*T)(\beta),\beta>\geq \Vert \beta\Vert^2,\]
so that
\begin{equation}\label{racine}
\Vert B(\alpha)\Vert \leq \Vert \alpha\Vert_B,\; \forall \alpha \in E.
\end{equation}

In the sequel we will denote by $\rm{BUC}(E_B,\R)$ the space of bounded and uniformly continuous functions on $E$ equipped with the norm $\Vert . \Vert_B$. Let us remark that because of \pref{controlhprime}, bounded functions which are uniformly continuous on $E$ equipped with the usual norm of $\R^d\times (H^1)^{\prime}$ belong to $\rm{BUC}(E_B,\R)$. In particular if $A$ is $H^1$ then the value function $v$ defined by \pref{defofv} belongs to  $\rm{BUC}(E_B,\R)$.

\smallskip

We end this paragraph by remarking that equation \pref{hjb} now can (at least formally) be rewritten as 
\begin{equation}\label{hjb2}
\lambda v(\alpha) +H(\alpha, \nabla_x v(\alpha))-x\cdot \nabla_x v(\alpha)+\<T^*(\alpha), Dv(\alpha)>=0, \alpha\in D(T^*).\end{equation}
As in M. Crandall and P.L. Lions \cite{cranlio5}, we will take advantage of this structure (where $\<T^*(\alpha), Dv(\alpha)>$ has to be understood  as $\<\alpha, T( Dv(\alpha))>$) by  imposing restrictions on test-functions (typically of the form $\Vert \alpha-\alpha_0\Vert^2_B$) rather than on $\alpha$.

\subsection{Comparison theorem}

The comparison result for \pref{hjb} then reads as
\begin{thm}\label{compa}
Assume that  {\bf{(H1)}}, {\bf{(H'2)}}  and  {\bf{(H3)}} hold. Let $v_1$ and $v_2$ be in $\rm{BUC}(E_B,\R) \cap C^0(E_w,\R)$) respectively a viscosity subsolution and a viscosity supersolution of \pref{hjb} on $\R^d\times L^2$ then $v_1\leq v_2$ on $\R^d\times L^2$. 
\end{thm}

\begin{proof}
Let us define $M:=\sup_{E} (v_1-v_2)$, $B:=(I+T^*T)^{-1}$  and
 \[\Vert \alpha \Vert_B^2:= \<B(\alpha),\alpha>,\; \forall \alpha=(x,z)\in E\] 
as before. For all $\eps>0$, $\delta>0$ and $\alpha_1:=(x_1,z_1)$, $\alpha_2:=(x_2,z_2)$, in $\R^d\times L^2$, let us set $\theta:=(\eps,\delta)$ and
 \[\Phi_\theta(\alpha_1,\alpha_2):=v_1(\alpha_1)-v_2(\alpha_2)-\frac{1}{2\eps} \Vert \alpha_1-\alpha_2\Vert_B^2-\frac{\delta}{2}\left(\Vert \alpha_1\Vert^2+\Vert \alpha_2\Vert^2 \right).\]
 Let us also define:
 \[M_{\eps,\delta}=M_{\theta}=\sup \left\{\Phi_\theta(\alpha_1,\alpha_2),\; (\alpha_1,\alpha_2)\in E\times E\right\}.\]
The weak continuity and boundedness properties of $v_1$, $v_2$ ensure that the supremum $M_{\theta}$ is attained at some points $\alpha_i^\theta=(x_i^\theta, z_i^\theta)$ for $i=1,2$. 
 
 \smallskip
 
 Let us set
 \[\Phi_1(\alpha):=\frac{1}{2\eps} \Vert \alpha-\alpha_2^\theta\Vert_B^2+\frac{\delta}{2} \Vert \alpha\Vert^2,\; \forall \alpha\in E,\]
 we then have:
 \[ÊD\Phi_1(\alpha)=\left(\frac{1}{\eps}B_1(\alpha-\alpha_2^\theta)+\delta x,  \frac{1}{\eps}B_2(\alpha-\alpha_2^\theta)+\delta z\right),\; \forall \alpha=(x,z)\in E.\]
Since $\alpha_1^\theta$ is a maximum of $v_1-\Phi_1$ on $E$ and since $v_1$ is a viscosity subsolution of \pref{hjb}, we get from definition \ref{defvisc}:
 \begin{equation}\label{inegssvisc}
 \lambda v_1(\alpha_1^{\theta})+H(\alpha_1^{\theta}, q_{\theta}+\delta x_1^\theta)+\mathrm{liminf}_{\alpha=(x,z) \in E_0\rightarrow \alpha_1^\theta} \<\frac{1}{\eps}B_2(\alpha-\alpha_2^\theta)+\delta z, \dot{z}>\leq 0,
 \end{equation}
 where 
 \begin{equation}\label{defqtheta}
 q_\theta:=\frac{1}{\eps} B_1(\alpha_1^\theta-\alpha_2^\theta).
 \end{equation}
 We then remark that
 \[\<z,\dot{z}>=-\frac{1}{2} \vert z(0)\vert ^2\rightarrow -\frac{1}{2}\vert x_1^\theta\vert^2 \mbox{ as } \alpha \in E_0\rightarrow \alpha_1^\theta.\]
 Next, we  write:
 \[\begin{split}
  \<B_2(\alpha-\alpha_2^\theta), \dot{z}>=\<B(\alpha-\alpha_2^\theta), T^*(\alpha)>- B_1(\alpha-\alpha_2^\theta)\cdot x\\
  =\<TB(\alpha-\alpha_2^\theta), \alpha>- B_1(\alpha-\alpha_2^\theta) \cdot x
   \end{split}\]
 and since $TB$ and $B_1$ are continuous, we get:
 \[ \<B_2(\alpha-\alpha_2^\theta), \dot{z}> \rightarrow \<TB(\alpha_1^\theta-\alpha_2^\theta), \alpha_1^\theta>-B_1(\alpha_1^\theta-\alpha_2^\theta)\cdot x_1^\theta  \mbox{ as } \alpha \in E_0\rightarrow \alpha_1^\theta.\]
 Hence \pref{inegssvisc} can be rewritten as:
  \begin{equation}\label{inegssvisc2}
 \lambda v_1(\alpha_1^{\theta})+H(\alpha_1^{\theta}, q_{\theta}+\delta x_1^\theta)-q_\theta\cdot x_1^\theta- \frac{\delta}{2} \vert x_1^\theta\vert^2 
+\frac{1}{\eps} \<TB(\alpha_1^\theta-\alpha_2^\theta), \alpha_1^\theta>  \leq 0.
 \end{equation}
 Using in a similar way the fact that $v_2$ is a viscosity supersolution, we arrive at:
 \begin{equation}\label{inegssvisc3}
 \lambda v_2(\alpha_2^{\theta})+H(\alpha_2^{\theta}, q_{\theta}-\delta x_2^\theta)-q_\theta\cdot x_2^\theta+ \frac{\delta}{2} \vert x_2^\theta\vert^2 
-\frac{1}{\eps} \<TB(\alpha_2^\theta-\alpha_1^\theta), \alpha_2^\theta>  \geq 0.
 \end{equation} 
Substracting \pref{inegssvisc2} and \pref{inegssvisc3}  then yields:
 \[\begin{split}
 \lambda( v_1(\alpha_1^{\theta})-v_2(\alpha_2^{\theta}))+ H(\alpha_1^\theta, q_\theta+\delta x_1^\theta)- H(\alpha_2^\theta, q_\theta-\delta x_2^\theta)\\
 +q_\theta \cdot (x_2^\theta-x_1^\theta)-\frac{\delta}{2} (\vert x_1^\theta \vert^2+\vert x_1^\theta \vert^2) 
 +\frac{1}{\eps}  \<TB(\alpha_1^\theta-\alpha_2^\theta), \alpha_1^\theta-\alpha_2^\theta> \leq 0. 
 \end{split}\]
 Now, thanks to \pref{tbpositif}, the last term is nonnegative, which gives:
 \begin{eqnarray}
 \lambda( v_1(\alpha_1^{\theta})-v_2(\alpha_2^{\theta}))+ H(\alpha_1^\theta, q_\theta+\delta x_1^\theta)- H(\alpha_2^\theta, q_\theta-\delta x_2^\theta)\nonumber\\
 +q_\theta \cdot (x_2^\theta-x_1^\theta)\leq \frac{\delta}{2} (\vert x_1^\theta \vert^2+\vert x_1^\theta \vert^2)\label{ineg4} . \end{eqnarray}
 By standard arguments (see for instance \cite{barles} or \cite{capuzzo}) and thanks to the fact that $v_1$ and $v_2$ belong to $\rm{BUC}(E_B,\R)$, one easily obtains:
 \begin{equation} \label{estimpert} 
 \lim_{\theta\rightarrow (0^+,0^+)} \left(\delta \Vert \alpha_i^\theta\Vert^2, \eps^{-1}\Vert \alpha_1^\theta -\alpha_2^\theta\Vert_B^2, M_\theta, v_1(\alpha_1^\theta)-v_2(\alpha_2^{\theta})   \right)=(0,0, M, M).
 \end{equation}
 On the one hand, using \pref{balal} and \pref{estimpert}, we have:
 \[ \vert x_1^\theta - x_2^\theta \vert \leq \sqrt{2} \Vert \alpha_1^\theta -\alpha_2^\theta \Vert_B=o(\sqrt{\eps})\]
 on the other hand, using \pref{racine}, we have
 \[ \vert q_\theta \vert \leq \frac{1}{\eps} \Vert B(\alpha_1^\theta-\alpha_2^\theta)\Vert \leq  \frac{1}{\eps}  \Vert \alpha_1^{\theta}-\alpha_2^{\theta} \Vert_B=o(\frac{1}{\sqrt{\eps}})\]
 so that 
 \[\lim_{ \theta\rightarrow (0^+,0^+)} q_\theta \cdot (x_2^\theta-x_1^\theta) =0.\]
 Lemma \ref{modulh} and  \pref{controlhprime} then imply that there is a nonnegative constant $C$ such that:
 \begin{equation}\label{ucdeH1}
 \vert H(\alpha,p)-H(\beta,p)\vert \leq C \Vert \alpha-\beta \Vert_B (1+\vert p\vert), \; \forall (\alpha,\beta, p)\in E^2\times \R^d,
  \end{equation}
 and 
  \begin{equation}\label{ucdeH2}
 \vert H(\alpha,p)-H(\alpha,q)\vert \leq C \vert p-q \vert  (1+ \Vert \alpha \Vert ), \; \forall (\alpha,p, q)\in E\times \R^d\times \R^d.
\end{equation}
 We thus deduce
\[\begin{split}
  H(\alpha_1^\theta, q_\theta+\delta x_1^\theta) &- H(\alpha_2^\theta, q_\theta-\delta x_2^\theta)=
 H(\alpha_1^\theta, q_\theta+\delta x_1^\theta)- H(\alpha_2^\theta, q_\theta+\delta x_1^\theta)\\
&+ H(\alpha_2^\theta, q_\theta+\delta x_1^\theta)- H(\alpha_2^\theta, q_\theta-\delta x_2^\theta)
\\
&\leq C \left(  \Vert \alpha_1^\theta-\alpha_2^\theta\Vert _B (1+ \vert q_\theta\vert +\delta (\vert x_1^\theta\vert +\vert x_2^\theta\vert)) \right)\\
&+ C \left( (1+ \Vert \alpha_2^\theta\Vert ) \delta  (\vert x_1^\theta\vert +\vert x_2^\theta\vert) \right)\\
&=o(\sqrt{\eps}) o(\frac{1}{\sqrt{\eps}})+o(\frac{1}{\sqrt{\delta}})o(\sqrt{\delta})\rightarrow 0 \mbox{ as } \theta\rightarrow (0^+,0^+) .
 \end{split}\]
Putting everything together and passing to the limit in \pref{ineg4} then yields $\lambda M\leq 0$ so that the proof is complete.

\end{proof}

We have already noticed that if $A$ is $H^1$ then $v$ defined by \pref{defofv} actually belongs to  $\rm{BUC}(E_B,\R)$. We thus deduce the following

\begin{thm}\label{mainthm}
Assume that  {\bf{(H1)}}, {\bf{(H'2)}}  and  {\bf{(H3)}} hold, then the value function $v$ defined by \pref{defofv} is the only $\rm{BUC}(E_B,\R) \cap C^0(E_w,\R)$ viscosity solution of \pref{hjb} on $\R^d\times L^2$.
\end{thm}

\section{Variants and concluding remarks}\label{sec5}

\subsection{Finite-dimensional reduction}

For the sake of simplicity, let us assume in this paragraph that $d=k=1$ and that $A$ is an exponential weight:
\begin{equation}\label{aexpo}
A(s)=e^{-\delta s}, \; \forall s>0 \mbox{ with } \delta>0.
 \end{equation}
Defining the optimal control problem and its value function $v$ as in \pref{defofv}, it is easy to see that, under the special exponential form of $A$, $v$ actually depends on $z$ only through the scalar parameter $y(z):=\int_0^\infty e^{-\delta s} z(s)ds$. More precisely, setting:
\begin{equation}\label{vw}
v(x,z)=w(x,y(z)), \; \forall (x,z)\in E \mbox{ and } y(z):=\int_0^{+\infty} e^{-\delta s} z(s)ds.
\end{equation}
it is easy to check that $v$ solves the infinite-dimensional Hamilton-Jacobi equation \pref{hjb} if and only if $v$ is given by \pref{vw} and $w$ solves the two-dimensional equation:
\begin{equation}\label{hjbw}
\lambda w(x,y)+H_0(x,y,\partial_x w(x,y))-\partial_y w(x,y) (\delta y+x)=0,
\end{equation}
where 
\[H_0(x,y,p):=\sup_{u\in K} \{-L(x,u)-p\cdot F(x,u,y)\},\; \forall (x,y,p)\in\R^3.\]
This finite-dimensional reduction of the problem of course heavily relies on the exponential form \pref{aexpo}. We refer to \cite{larssen} for the extension of such finite-dimensional reduction in a stochastic setting .

\subsection{The evolutionary problem}

In the present article, we have focused on the stationary case. If we consider, the finite horizon optimal control problem
\begin{equation}\label{defofvt}
v(t,x,z):=\inf_{u\in V} \left\{ \int_t^{T}  L(s,y_{t,x,z,u}(s),u(s))ds+g(y_{t,x,z,u}(T))\right\}
\end{equation}
where $y_{t,x,z,u}$ denotes the solution of the Cauchy problem
\[\begin{split}
\dot{x}(s)=F\left(s,x(s),u(s), \int_0^{+\infty} A(\tau) x(s-\tau) d\tau\right), \; t>0,\\
\mbox{ with initial conditions }  x(t)=x, \; x(t-s)=z(s), s>0,\end{split}\]
this leads to the following evolution equation for $v$:
\[\partial_t v(t,x,z)+\inf_{u\in K} \left\{L(t,x,u)+\nabla_x v(t,x,z)\cdot \F(t,x,z,u)  \right\}-\<D_z v(t,x,z),\dot{z}>  =0,\]
together with the boundary condition
\[v(T,x,z) = g(x),\; \forall (x,z)\in E.\]

\end{document}